# New Lower Bounds for Van der Waerden Numbers

Alexey V. Komkov

October 24, 2020

**Abstract**

This work contains certificates numbers Van der Waerden, was found using SAT Solver and Genetic Algorithm. These certificates establish the best currently known lower bounds of the numbers Van der Waerden W( 7, 3 ), W( 8, 3 ), W( 10, 3), W( 11, 3 ), W( 17, 3 ), W( 7, 4 ), W( 6, 5 ), W( 5, 6 ).

## 1 Introduction

Van der Waerden's theorem states that for any positive integers r and k there exists a positive integer N such that if the integers {1, 2, ..., N} are colored, each with one of r different colors, then there are at least k integers in arithmetic progression all of the same color. The smallest such N is the van der Waerden number W(r, k). Van der Waerden Numbers are quite difficult to calculate, at the time of writing the article, the exact values are known only for 7 van der Waerden numbers, and for the rest of the numbers only bounds are known. To prove the lower bound of a van der Waerden number it is sufficient to get a certificate of this number, i. e., a sequence of numbers {1, 2, ..., N} of r colors, and not having same-coloured arithmetic progressions of length k . Then N+1 will be the lower bound of W(r, k) number.

## 2 Certificates

Certificate W( 7, 3 ) > 343

```
00113663354603601544124523551354321606323600502552610664253321304402422360632310654614
46045015541236235466336112216245045262355052341346211242050110356556042446030452450310
01136633546036015441245235513543216063236005025526106642533213044024223606323106546144
6045015541236235466336112216245045262355052341346211242050110356556042446030452450311
```

Certificate W( 8, 3 ) > 515

```
25216145140707307525447527532400405716505705166436422745412427346146762526766700547542
36574234243550556032256056611542531076533534772672013232011050274273716073773546256210
73321121163653674001674677403203064373063025307304461004404651321365474365363271671455
06145140707307525447527532400405716505705166436422745422427346146762526766700547542365
74234243550556032251056612542533076533534772672013231012050274273711073773546251210733
2112126365367420167467740320306437306302560730446100440465132136547436536327167141201
```

Certificate W( 10, 3 ) > 892

14318463515917437942835565522921902821127412947411511798800422442893511853372711892701
34064040150389383094765092984150786832352533162893093771274060223554175252296349064350
15976978656678774426608428929355459476606657524687366328426338308104254037223634608184
47760855861937707503497943898902630405688550438697691727910795905917861511466920611420
66566350054041141153866153801366060336844294756533760174142254098196852199893580376249
07473236038425187084635159174379428355655229219028211274029474115917958004224428935118
53772711812701740640401503813830947650929841507867723525731627930937712740602235541752
52296349064350159769786566787744266084289293554594766066575246873663284263383081042540
37223684608184473608558619377075034179438989026304057885504386976117279107959059178615
01466920691420665663500540411411538661538013660618768442947565337601741422540951968521
99890880376249079772360384259530

Certificate W( 11, 3 ) > 1187

48466199AA74A6788289168496059204409AA4A563091073973252A14A11951012A1758A297550982913A8
A51A83257A89224125488355188058AA7095A91256428A560A50037930A90A52273973034077582A840970
5451794A63485A867A93A4493764A03028159106265223A690013247097281883077237570574851991930
037951963684263810650971683982A53582254A34892673832643681510880953094683011817A60A4121
1854570577A9756365169726789396028AA92882061271690A6A787A1463323705A844A8A2179748672646
32677544A63513592830230A17A6083007A69407291A81054A670651873654891A9225A572A5A97A633912
61087729A64793296101193493289351356126307A83A26406227001468125265703453A21108866199AA7
4A6788289168496059206409AA5A563091073973252A14A11961012A1758A297550982913A8A51A83257A8
9223125488355188058AA7095A94256428A560A50037930A90A52273973034077582A8409705451799A634
85A867A93A4493764A03028159106264223A69004325709728188307723757057485199193003795196368
4263810640971683982A53582254A44892673832613681510880953094683011817A60A11211844570577A
9756365169726784396028AA92882061271690A6A187A1463323705A844A8A247974867264632677564A63
513592840230A17A4083007A69407291A84053A670651873654891A9225A572A5A97A63394261087724A67
793296101193493289351356126304A83A26406227001468125265703454A21104A44

Certificate W( 17, 3 ) > 3549

3AG0173AG38FC1BDB0E6G0G031A1652107AG3375157ABC827952G09CFC6731BCC5C6B4BG46D58E3B00C2CG
293GE95592DA049G180DBD13GDGE9BG9847823301BF1G241A3342EB3FE364C862CD0822EA7C61E5CG07A6C
3F56B3AG89FD82D16362CG03F68572EF2BB16G2165G04970FF73D43D2FB2BBDE229B2F618459ADGB06G7EE
7F603G7155B44B296F5GE08G0G53B09325G62G499FDF65CDE950655BB1C49GF3GDC6F0G573BECDC1E83B4C
GD8A4CE18DCFB91D83G0FB216D6271F2BG965DG671B48G0FF81G2EDC605FD49D282DG49E716B3C4DBD1105
F088D51A81D1C25F6DFA494C6B81AE942925A4C85762D3E19709E5GBE82080DAB9651829C09377925FD0C5
06EDC646709B3F00BG2673F40959F593D1BG1157CGD971107104842605BB43508B7F41971677B98BG7F28C
43BD9G6DAB8AF71F9A0CBF2F28GG1B22F77CD0821868343CG7177303CD5G22A61347G633D9E7562DF01879
F62FB43F9G54GBEA45G28CE0377664144A751F8GAC9D7EE8AA9351953CDC9F29B9B03EF2FB7G41CD7F25F4
E3791B8AD3EF7EB50GC941G74ABCB50862FG8642652G6C6D1G7DBDC8D1EB4CE18F9873AB0942EB3GB0B20G
A21DDC408E3CB8A27G68GD1E518A8F4E65095DF08C20G793EDA47857A48A2267D363B4C6723BE3F59E4C63
6D5155GB86D11BDAD0B41G28A0EA52172CB864C120B69CF0E52F6C37D808BG4CD79AD164BF0E65C61BFC15
43D165DEF2FE9CDC1FDFAC870F07C00C796E9FG884087DCD83FDF78A04375968GG0870E53308C499B2AG88
D7E631BCCG7EA3F0D6256541B4B2292DFG4BF2162CACCF0G72G083D2EA4FG4CG2G151B067C9A816DA5DF75
B40B45F6714C67C8291DC6G2953CD435C6A83B0BG2668361A71FB4GE0A57658BG1C7DE64CGD651E72G4DE2
587DBE1A7C9E3GG0C55102BG1B4C8F7F6094CFF8G04827957F4C8D5G68B4C6D28A9360670E0DGGDBDCF58B
409D0727FA4EG5F77108351E25AACC218B7A64EG9G57C649D59E9C4D3G7F73FEC9D78731D8AF115D3A8347
975BD0A5C859B043FCEADA0882DE653D8EG56EDBA0322FA5BDCCGD6B81E9CG167DG418459B8A61B20875A6
4BD9F531A4A7CF6DB0CBG1975AA9B74GB8348921EB7GGE3FCG223608F91B21489509B9143B25654A0E941E
E8363B4DEF7CG86F29633GDC45E25CE83GD955G4E6GCF62GG185C823CG89276653ED8D1F25C609571773B0
9F7E6DAD1C251721DC879E2214CB519F7348EG1865G6C9F5E223BCEC0173AG38FC1BDB0E6G0G031A165210
7AG3375157ABC827952G09CFC6731BCC5C6B4BG46D58E3B00C2CG293GE95592DA049G180DBD13GDGE9BG98
47823301BF1G241A3342EB3FE364C862CD0822EA7C61E5CG07A6C3F56B3AG89FD82D16362CG03F68572EF2
BB16G2165G04970FF73D43D2FB2BBDE229B2F618459ADGB06G7EE7F603G7155B44B296F5GE08G0G53B0932
5G62G499FDF65CDE950655BB1C49GF3GDC6F0G573BECDC1E83B4CGD8A4CE18DCFB91D83G0FB216D6271F2B
G965DG671B48G0FF81G2EDC605FD49D282DG49E716B3C4DBD1105F088D51A81D1C25F6DFA494C6B81AE942
925A4C85762D3E19709E5GBE82080DAB9651829C09377925FD0C506EDC646709B3F00BG2673F40959F593D
1BG1157CGD971107104842605BB43508B7F41971677B98BG7F28C43BD9G6DAB8AF71F9A0CBF2F28GG1B22F
77CD0821868343CG7177303CD5G22A61347G633D9E7562DF01879F62FB43F9G54GBEA45G28CE0377664144
A751F8GAC9D7EE8AA9351953CDC9F29B9B03EF2FB7G41CD7F25F4E3791B8AD3EF7EB50GC941G74ABCB5086
2FG8642652G6C6D1G7DBDC8D1EB4CE18F9873AB0942EB3GB0B20GA21DDC408E3CB8A27G68GD1E518A8F4E6
5095DF08C20G793EDA47857A48A2267D363B4C6723BE3F59E4C636D5155GB86D11BDAD0B41G28A0EA52172
CB864C120B69CF0E52F6C37D808BG4CD79AD164BF0E65C61BFC1543D165DEF2FE9CDC1FDFAC870F07C00C7
96E9FG884087DCD83FDF78A04375968GG0870E53308C499B2AG88D7E631BCCG7EA3F0D6256541B4B2292DF
G4BF2162CACCF0G72G083D2EA4FG4CG2G151B067C9A816DA5DF75B40B45F6714C67C8291DC6G2953CD435C
6A83B0BG2668361A71FB4GE0A57658BG1C7DE64CGD651E72G4DE2587DBE1A7C9E3GG0C55102BG1B4C8F7F6
094CFF8G04827957F4C8D5G68B4C6D28A9360670E0DGGDBDCF58B409D0727FA4EG5F77108351E25AACC218
B7A64EG9G57C649D59E9C4D3G7F73FEC9D78731D8AF115D3A8347975BD0A5C859B043FCEADA0882DE653D8
EG56EDBA0322FA5BDCCGD6B81E9CG167DG418459B8A61B20875A64BD9F531A4A7CF6DB0CBG1975AA9B74GB
8348921EB7GGE3FCG223608F91B21489509B9143B25654A0E941EE8363B4DEF7CG86F29633GDC45E25CE83
GD955G4E6GCF62GG185C823CG89276653ED8D1F25C609571773B09F7E6DAD1C251721DC879E2214CB519F7
348EG1865G6C9F5E223BCEC

Certificate W( 7, 4 ) > 9980

010023455031053046544322431124204353405254461145420061213644130466546430425015325441263164045265
065463445315446440143260424455455224325513035224064335416062516162020344100164653216251244044223
214424401665553153255202663556355503250121131331312026315551516212023653046335005153223626020615
013146330166322535524642664345344610066060052202545113636052426265602114016021410630351130121531
134603260655113240626332603042063200260212456416260454141261232166306112040015500256200021403261
243210334454161223315640110600631460660561435605405640103366131054513620255255215623206401024100
666160313116421053224503310145024564522203313462233652210434103121441024146033321126320030621231
331615562441054021566463303613005024440153151205432411266101235010022121005023423404636251423443
024024636350306031024332051402325256162142303103224220233664602520006553032324200130424110065136
041214003323131460636231304633030312163662546626333022522326501230625621565326654131641262512661
612546466024302511465043361030226266623616520123046361056160033565560116012351442433240106312043
564140304323113010354142630135226212060043612014133020336263230066614261600106022614113511135606
031013526303023420661163332013330331606233230366315156261012160332300350130622251321400133001630
161433513161644301513531266262200630615040444116122544543616202063220366626623336012306333012441
601201063334005056302231021462223523210003360235212365230011200311412342431550036035410414646312
200261230650511055126323132255241134230241101401216036152504140363006054315605306106132614302314
010062040621061312544666522602620023202154636314201506110151042420331422500606343426254261361156
246231053436222115160113311535465622531556051132162562561013221120153602334212555202405343430452
461616250206660306110241632033001323440560430133425663160056564163311332534330004101623105504633
513506604051614252240064112000623632334263206315145630201022050250065456560313523052620422320310
216350134235216246315624221412645240635124303253233502104133606600131014036423362642103402200315
645362356124122663140610343312203322106510304102511222301061000553023100214626226111523203253560
321114210325111055045512001520205351463460035336452626352152200505503161363216443431341164532153
221132243013600052132612210112530323505636315412101105253112111320111433442051021214061323125244
133361336340024232243506355522101505112021136320451422520305006132140631612321331016212636541020
315232601110663610325332566461122534423515621032045351055505002123511626533360216024465560343550
360403643136655016563046052103265010040021115055560551530312121156213104142554131011223303625136
422330052132205010121464222060246451102006001231210630435604010213620116023425212615154026021641
063604145261021026223030022326103235503300146203631632646026221351211564443026623350566343113103
052122015033011124536023643031236162300641100511341200046566026632311265001623064331312535552230
023255201046306606402051364623135511232656266246163542452563152252332656311003533666112301103501
263022205602266322420335214530103503635662505355660302502401520621240114054501336652501251133163
032133612152313025326330244505030625151336660200422525522355116314550446066160015542115531326352
025041001634226115021651020305363666315020313113133032601266615661550206601631646553236036341004
324404605341655544234412020544360445344412522006421346013004456663004401442565146316641526406506
325531503656016326424155465522132441343422044044564646244416202033110416664321424126601042424514
410166555315325520266355634550325012113133131202631545151621202365302633500515322362602061401316
633016632253452566266631533641004606004220255511363605212626560211501602131063035113012143113360
326065411324062633260305204320026021255661626055615126123216630411206001550025620002150326123321
033455116122331565011060043156066046143540550565010336613105551362025524521462320660102210066616
031311652104322440331414502556552220331336224364221053510312155102516603332112632003062123133461
556261105402156646330361300402545015315120443221126610123401002212100502302340663625142345302402
663435030603102433205130232524616215430310322522023346560244000655303232620013042211006513604121
100332313166043623130563303031216366244642633302252232650123062562156432465513164124251266161246
646602530251116504336103022626662361452012405636105616003356456011601235155233325010631202356514
030532411301035315263013422621204005361201513302033626323006661526140010602261511351113560603101
352630302362066116333201333033160623323036631515626101216033230035013062225132450013300163016153
351316165530151343126626220043061505066511612251156361620206322036662662333601230633301246160120
106333500505630223102166222352321000336023521236523001120031161235253154003603561051464631220026
123064041105512632313225525113524026110160121603615240615036300405431560530610613261530231201006
205062106131256466452260262002320215663631520150611015105222033152250060635362625526136115426623
105323622211516011331153556562253154605113216256256101322112015360233521255520260535353065216161
625020666030611025163203300132356056053013362566316005656616331133253633000510162310550563351350
660405161525226006311200062363233424320631515563020102205025006555656031352305262052232031021435
013523521625631562522141263526063512530325323340210613360660013101503662334245210350220031565536
235612412266316041034331220332210651030610251122230106400055302310021542622611142320325346032111
521032511105505451200152020435156336003433645262635215220040550314136321665343134116343215322113
226301360005213261221011253032350543431551410110525411211132011153344205102121504132312525613336
133634002523225350435554210150511202113632035162242030500613215063161232133101621263654102031423
260111066361032433256646112253462341562103204535405550500212351142653334021602546556035344036050
365313665501656305605210326401006002111504456054153041212115621310515245313101122330362413652233
005213220501012166522206024654110200600123121063053560601021362011602312424261414402602166106360
515526102102622303002232610323550330016620363163266602622135121156556302662335056634311310305212
201504301112553402366303123616230064110051131120004656602663231426500162306433131254545223002324
420104630660660205136642313551123265626625616356245246315225233245631100353366611230110350126302
220540226632252033421453010350363566250435564030250420142062125011505250133665250125113316303213
361214231302532633025650503062414133666020062242552435511431444045606616001453211553132635202505
100163622611502165102030536366631502031311313303260126661566155020660163166654323603630444232554
360530144453223314204453546452435611252440612136601440665666300250453256542441460144450650632543
440365644634602215556552413255440454400603346160645164640443311004466544464512664102232151241014
455531532552026635563555032501211313313120263154515162120236530243350051432236260206150131463301

```
663224344246626664153364140660600522024551136360521262656021150160214106303411301215311336032606
551132606263326030420632002602125566162605561512612321663061120600155002542000215032412332103355
511612233156401106006315606605615356055056501033661310554136202552452144232066010221004661603131
164210532264033101540245454222033133422336424105351031215510251660333211263200306212313316154626
110540215645633034130050254404431412054322112461012350100221210050230235066342414234430250246363
503060310253320513023252561621523031034252202334656024200065530323262001306221100651360512110033
231316606362313056330303121646625666263330224223265012306256214653264541316412625126614125664660
253025111640533610302262666236165201230563610561600334455601160123515523332501063120235651503053
231130103431526301352262120600536120451330203362432300666142616001060224141135111356060310135263
030236206611633320133303314062332303663151562610121603323003401306222514215001330016301615334131
616453015134312662622006306150504651161225115636162020634203666266233360123063330125616012010633
350050563022310216622234232100033602352123654300112003114123524314400340346105166563122002612306
505110541263231422552414352302611016012160361525061503630060463156053061061424153023120100620406
210613125656665226026200232024466363154015061101510522203315224004063536242544413411562462310432
362221451401133114355656225315560511321625624610132211201546023352125552026054534306521614162502
066603061102516320330013435605604301336254631600565461633113325463300051016231045046334135066050
516142522400641140000624632334263206315155630201022040250064556560413523052620522320310216350435 2
352162563156252214126342606341253032432334041061336066001310150366233624521035022003156453623561
241246631606103434122033221064103061025112223010610005530231002156262241115234032535603211152103
251110450555140015202053515634600343364426263521422005045034613632146535313611634321532211322630
136000521326122101124304235056344155121011052531121113201115335520510212150413231254541333613363
500244322435063454221015051120211363203416224203050041321404316123213310142126364610203142326011
106634103253325664611224356235156210320653510555050021235114265333602160255645603535503404036431
366540164430460521032640100600211150555604515303121241562141041424531310112243036251364223300521
322040101246652220602464511020060012312106304356040102134201160231242126151540260214610636041552
610210262230300223261032354033001662036314326660262213512115654640264233504464531431030421220150
330111244360236630312461623006511004113112000664660266323112650016230453313125354522300232552010
463066066020513666231355112326462662561645425425631422523326563110035336461123011034012630222056
022663224203352145301035036346624043456404025022015206212501150525013366524012411331630321336121
423140253263304564050406251443344602006224244224551163145405540664400155321155313263520250510016
362461150216410203043636463150203131131330324012666156615502066016316665532360363010024255036444
016555324331140205353605253561424220441444660430066464644025415344641263466015265144044255415034
56016424022154564442144551303522006033561404251614242043144016645321625456601025341512514466
```

## 3 Results *

| Number W( r, k ) | Old Lower Bounds | New Lower Bounds |
|---|---|---|
| W( 7, 3 ) | > 342 [1] | > 343 |
| W( 8, 3 ) | > 511 [1] | > 515 |
| W( 10, 3 ) | > 889 [1] | > 892 |
| W( 11, 3 ) | > 1183 [1] | > 1187 |
| W( 17, 3 ) | > 3546 [1] | > 3549 |
| W( 7, 4 ) | > 9951 [2] | > 9980 |
| W( 6, 5 ) | > 99483 [2] | > 99554 |
| W( 5, 6 ) | > 540170 [2] | > 540197 |

* Certificate W( 6, 5 ) and W( 5, 6 ) in http://komkov.org/VanDerWaerden/